\documentclass[final,3p,times]{elsarticle}

\usepackage{amssymb}
\usepackage{amsmath}

\usepackage{hyperref}

\usepackage{booktabs}
\usepackage{subcaption}
\usepackage{bm} 
\usepackage{upgreek}
\usepackage{siunitx} 
\usepackage{empheq} 

\usepackage{notations}

\usepackage{lineno}

\journal{Computer Methods in Applied Mechanics and Engineering}

\begin{document}

\begin{frontmatter}



\title{A space-time LATIN-PGD strategy for solving Newtonian compressible flows}

\author[label1]{É. Foulatier\corref{cor1}} 
\author[label1]{P.-A. Boucard}
\author[label1]{F. Louf}
\author[label1]{D. Néron}
\author[label2]{P. Junker}

\cortext[cor1]{Corresponding author}

\affiliation[label1]{organization={Université Paris-Saclay, CentraleSupélec, ENS Paris-Saclay, CNRS, LMPS - Laboratoire de Mécanique Paris-Saclay},
            city={Gif-sur-Yvette},
            postcode={91190}, 
            country={France}}

\affiliation[label2]{organization={Leibniz University Hanover, Institute of Continuum Mechanics},
            city={Hanover},
            country={Germany}}

\begin{abstract}
Simulating flow problems is at the core of many engineering applications but often requires high computational effort, especially when dealing with complex models. This work presents a novel approach for resolving flow problems using the LATIN-PGD solver. In this contribution, we place ourselves within the framework of Newtonian compressible and laminar flows. This specific and relatively simple case enables focusing on flows for which a state equation provides a direct relation between pressure and density. It is then possible to use the LATIN solver to set up a pressure-velocity decoupling algorithm. Moreover, Proper Generalised Decomposition (PGD) is natively included in the solver and yields two independent space-time decompositions for the velocity and the pressure fields. As a first step, the solver is validated on a problem for which an analytical solution is available. It is then applied to slightly more complex problems. The results show good agreement with the literature, and we expect that the solver could be used to compute more complicated material laws in the future.
\end{abstract}



\end{frontmatter}


\section{Introduction}\label{sec:intro}

Despite the formalisation of Navier-Stokes equations in the first half of the 18th century, analytical solutions were only obtainable for particular cases for many years. During the 1950s, advancements in the field of computational science enabled the resolution of the intricacies inherent within these equations. The first computational approach to solving incompressible flow problems utilised finite differences \cite{chorin_1967, chorin_1968} while Taylor \cite{taylor_1973} proposed a method that employed finite elements a few years later, thereby allowing the management of more complex geometries.

Numerical solvers dedicated to flow problems are constantly evolving to meet specific needs. For example, recent studies have employed various methods, including advanced stabilisation techniques \cite{arndt_2015}, log-conformation reformulation for viscoelastic flows \cite{habla_2014}, and the PFEM for evolving domains \cite{cremonesi_2020}. Direct coupling methods, also known as pressure-velocity formulations, are widely used to solve fluid problems. These methods involve gathering all the equations of the problem into a single matrix, which can lead to high computational costs, as mentioned in \cite{wang_2018}. Therefore, one can resort to pressure-velocity decoupling strategies presented in \cite{haroutunian_1993,benim_1986}.

However, all these approaches require high computational effort, mainly due to the problem's non-linearities. In this context, it is advisable to investigate surrogate models or model order reduction techniques. Despite being non-intrusive, surrogate approaches (such as metamodels or neural networks) are designed to compute specific quantities of interest, and one needs to fully reconstruct the model if they want to access other data. For in-depth studies where we compute many fields, building a reduced-order model seems more appropriate. In most cases, model order reduction techniques consist of two stages: the \textit{offline phase}, which is computationally expensive, and the \textit{online phase}, which is significantly faster. The offline phase involves computing several high-fidelity solutions, referred to as \textit{snapshots}, for various parameter sets. One uses the snapshots to compute a reduced-order basis, and once it is obtained, the online phase enables the rapid computation of solutions for new problems. \textit{Proper Orthogonalised Decomposition} (POD) \cite{chatterjee_2000} is a widespread model order reduction technique. During the online phase, it involves projecting the equations of the new problem onto the reduced basis to obtain a fast and accurate solution. We encounter POD in many CFD problems \cite{berkooz_2003,mendez_2019,beckermann_2025}, but methods like gappy-POD \cite{willcox_2006}, reduced basis \cite{quarteroni_2007}, or manifold interpolation \cite{franz_2014} might also be suitable. All of these methods prove to be efficient in parametric studies, but the offline phase becomes costly and not worthwhile in the frame of studies with few computations. Therefore, one can imagine to investigate model order reduction techniques, which provide a progressive enrichment of the solution throughout the computation.

The \textit{Proper Generalised Decomposition} (PGD) \cite{ladeveze_1999,chinesta_2011} offers such a possibility in the frame of flow problems. In the first application of PGD for Navier-Stokes equations \cite{dumon_2011}, the separated variable decomposition focuses only on the space variables. Later, PGD has been applied to simulate harbour agitation \cite{modesto_2015} and the solution field is decomposed in space functions and functions depending only on the wave direction and frequency. The fluid resolution schemes are indeed generally incremental, and it is therefore impossible to build a solution across the entire space-time domain without considering all time steps. This issue is crucial for the application of PGD, and one must develop non-incremental schemes to allow a space-time decomposition. A first attempt is presented in \cite{aghighi_2013}, where a fixed-point algorithm is used to generate new modes. Performing a linearisation around the solution obtained at the previous iteration enables us to consider non-linearities. However, the authors recognise that this linearisation technique is not optimal \cite{ammar_2010}. In this context, looking for an appropriate non-incremental and non-linear solver seems essential.

In the 1980s, Pierre Ladevèze proposed a novel approach, known as the LATIN method, to address non-linear structural mechanics problems \cite{ladeveze_1985, ladeveze_1999}. Moreover, the method includes natively PGD and thus enables a space-time decomposition of the solution. The LATIN-PGD solver has been applied to various fields, including geometric and material non-linearities, contact, and composite damage, among others. In the domain of fluid mechanics, a strategy based on the LATIN solver combined with multiscale domain decomposition methods has been developed for fluid-structure interaction \cite{vergnault_2010}. However, this approach does not make use of PGD. One can find more details on the domains of application of the LATIN solver in \cite{scanff_2021}. In particular, the LATIN-PGD version for multiphysics problems may be of interest here, as coupling strategies between velocity and pressure are crucial in the Navier-Stokes equations. Initially presented in \cite{dureisseix_2003} for poroelasticity problems, the solver relies on the separation of the difficulties. Each iteration of the algorithm consists of two stages: one that utilises only the coupled non-linear constitutive relations and another that employs global admissibility equations for each physical fields independently. The method proves to be very modular, enabling the construction of an independent basis for each primal fields and the use of different discretisations in time and space \cite{neron_2007,neron_2008}. More recently, the solver has been extended to more complex multiphysics problems (thermo-poroelasticity) in a more robust version \cite{foulatier_2025}. This latter turns out to be favourable when tackling parametrised problems.

This paper presents a novel approach based on the LATIN-PGD method to solve the Navier-Stokes equations. To the best of the authors' knowledge, such work has not been previously presented in the literature. It could overcome several difficulties related to non-linear flows: reducing the complexity of the problem through a decoupling strategy, enhancing computational efficiency with PGD, dealing with non-linearities, and introducing complex material laws at the local stage. To demonstrate the solver's capability to solve flow problems, we focus here on relatively simple 2D and 3D problems, specifically Newtonian laminar flows. Moreover, to have a constitutive law that couples pressure and velocity, the study focuses on compressible flows.

The paper is organised as follows. In \autoref{sec:ref_pb}, the continuous governing equations of the reference flow problem are presented, followed by their numerical discretisations in \autoref{sec:discretisation}. \autoref{sec:latin-pgd} focuses on the LATIN-PGD solver in the frame of fluid studies. In \autoref{sec:results}, we prove the relevance of the method on three numerical examples. The first example consists of a simple case for which an analytical solution is available. We then reproduce the well-known benchmark of a flow around a cylinder. Finally, one shows the extendability of the code to 3D geometries. \autoref{sec:conclusion} yields concluding remarks and possible outlooks for this work.
\section{Reference flow problem}\label{sec:ref_pb}

This section describes the reference problem by presenting the equations related to compressible Newtonian flows. We carry out the study on a bounded domain $\dom \subset \mathbb{R}^3$ during the time interval $I = [0, T]$. Equations of flows are derived from conservation principles \cite{schlichting_2016}. The mass conservation yields the continuity equation \eqref{eq:continuity}:
\begin{equation}\label{eq:continuity}
 \derivPart{\p}{t} + \divergVec{(\p\V)} = 0 
\end{equation}

\noindent where $\p$ is the fluid density and $\V$ the flow velocity.

The second flow equation comes from the conservation of momentum \eqref{eq:momentum}:
\begin{equation}\label{eq:momentum} 
 \divergTens{\Sig} + \fd = \p\derivFlu{\V}{t}
\end{equation}
with $\Sig$ the symmetric stress tensor, $\fd$ the volumic body forces. $\derivFlu{\V}{t}$ designates the material derivative of the velocity and is defined as \eqref{eq:mat_derivate}:
\begin{equation}\label{eq:mat_derivate}
 \derivFlu{\V}{t} = \derivPart{\V}{t} + \gradVec{(\V)} \V 
\end{equation}

Constitutive relations close the system of equations. The expression for the stress tensor is given in \eqref{eq:stress_law} where $\bm{\tau}$ is the deviatoric stress and $p$ is the pressure.
\begin{equation}\label{eq:stress_law}
 \Sig = \bm{\tau} - p\Id
\end{equation}
As emphasised in the introduction, the study is restricted to Newtonian compressible flows. In such case, the deviatoric stress is given by \eqref{eq:deviatoric_stress}, where $\Viscous$ is the viscosity tensor of order four and $\mu$ and $\lambda$ respectively designate the dynamic viscosity and the second viscosity \cite{rosenhead_1954}:
\begin{equation}\label{eq:deviatoric_stress}
 \bm{\tau} = \Viscous \mathbin{:} \GradSymV = 2\mu\GradSymV + \lambda \Tr \GradSymV \Id
\end{equation}

The second constitutive equation relates to the density $\p$. For compressible flows, a state equation links density and pressure. Here, we use the ideal gas equation \eqref{eq:ideal_gas_law} where $R$ is the ideal gas constant and $M$ the molar mass:
\begin{align}\label{eq:ideal_gas_law}
     p = f(\p) = \frac{R}{M} T \p = r T \p
\end{align}

For air, one can take $r = \qty{287.04}{\joule\per\kilogram\per\kelvin}$ \cite{oosthuizen_2013}. Again, we assume that the temperature is constant, and set $T = T_0 = \qty{293}{K}$.

Finally, one needs boundary and initial conditions. Boundary conditions are of Dirichlet type \eqref{eq:BC_Dirichlet_p}, \eqref{eq:BC_Dirichlet_v} and Neumann type \eqref{eq:BC_Neumann}:
\begin{subequations}
\begin{empheq}[left=\empheqlbrace]{align}    
    \label{eq:BC_Dirichlet_p} &p = p_d  \quad \text{over}\: \bord{p}\times\I \\
    \label{eq:BC_Dirichlet_v} &\V = \V_d  \quad \text{over}\: \bord{v}\times\I \\
     \label{eq:BC_Neumann} &\Sig \normal  = \Fd = -p_d \normal \quad \text{over}\: \bord{F}\times\I
\end{empheq}
\end{subequations}
The initial conditions \eqref{eq:IC_density} and \eqref{eq:IC_velocity} are for the density and the velocity, since the flow is compressible:
\begin{subequations}
\begin{empheq}[left=\empheqlbrace]{align}
     \label{eq:IC_density} & \p(t=0, \vect{x}) = \po = f^{-1}(p_0) & \forall \vect{x} \in \dom \\
     \label{eq:IC_velocity} & \V(t=0, \vect{x}) = \vect{0} & \forall \vect{x} \in \dom
\end{empheq}
\end{subequations}
where $p_0$ is the initial pressure in the medium.

\section{Numerical discretisation}\label{sec:discretisation}

This section deals with the numerical discretisation: we first discretise the space and then the time. Here, $\spaceH$ designates the Sobolev space of square-integrable functions whose first-order derivatives are also square-integrable.

Let us define the following function spaces:
\begin{itemize}
    \item $\spaceV = \{\, \V \mid \V \in \spaceH, \V = \Vd \: \text{over} \: \bord{v} \,\}$;
    \item $\spaceVo$, the associated homogeneous space ;
    \item $\spaceP = \{\, \p \mid \p \in \spaceH, \p = \pd \: \text{over} \: \bord{\rho} \,\}$;
    \item $\spacePo$, the associated homogeneous space.
\end{itemize}

The balance equations \eqref{eq:continuity} and \eqref{eq:momentum} given in \autoref{sec:ref_pb} can be written in a weak form using a variational formulation. Therefore, after integration, the compressible flow problem is:

Find $\V \in \spaceV$, $\p \in \spaceP$ such that $\forall t \in \I$: 
\begin{subequations}\label{eq:variational_formulation}
\begin{equation}\label{eq:var_density}
\displaystyle\intO{ \derivPart{\p}{t} \Ptest} + \intO{\scalprod{\rho\V}{\grad{\Ptest}}} = 0 \qquad \forall \Ptest \in \spacePo
\end{equation}
\begin{equation}\label{eq:var_velocity}
\displaystyle\intO{\rho \derivFlu{\V}{t}\Vtest} + \intO{\Sig : \GradSymV^\star} = \intO{\scalprod{\fd}{\Vtest}} + \int_{\bord{F}} \scalprod{\Fd}{\Vtest} \dS \qquad \forall \Vtest \in \spaceVo
\end{equation}
\begin{equation}\label{eq:var_vo}
 \V(\M, t=0) = \V_0 \quad \forall \M \in \dom 
\end{equation}
\begin{equation}\label{eq:var_po}
 \p(\M, t=0) = \po \quad \forall \M \in \dom 
\end{equation}
\end{subequations}

When replacing the stress tensor and the pressure by their expressions given by the constitutive relations \eqref{eq:stress_law} and \eqref{eq:ideal_gas_law}, one gets the subsequent weak form of the problem:

Find $\V \in \spaceV$, $\p \in \spaceP$ such that $\forall t \in \I$: 
\begin{subequations}\label{eq:variational_formulation_dvlp}
\begin{equation}\label{eq:var_density_dvlp}
\displaystyle\intO{ \derivPart{\p}{t} \Ptest} + \intO{\scalprod{\rho\V}{\grad{\Ptest}}} = 0 \qquad \forall \Ptest \in \spacePo
\end{equation}
\begin{equation}\label{eq:var_velocity_dvlp}
    \begin{split}
        \displaystyle\intO{\rho \derivFlu{\V}{t}\Vtest} + & \intO{ 2\mu\GradSymV:\GradSymV^\star}+ \intO{\lambda \Tr \GradSymV \:\Tr{\GradSymV^\star}} \\
        & - \intO{f(\p) \Tr{\GradSymV^\star}} = \intO{\scalprod{\fd}{\Vtest}} + \int_{\bord{F}} \scalprod{\Fd}{\Vtest} \dS \qquad \forall \Vtest \in \spaceVo
    \end{split}
\end{equation}
\begin{equation}\label{eq:var_vo_dvlp}
 \V(\M, t=0) = \V_0 \quad \forall \M \in \dom 
\end{equation}
\begin{equation}\label{eq:var_po_dvlp}
 \p(\M, t=0) = \po \quad \forall \M \in \dom 
\end{equation}
\end{subequations}

The system to solve is a non-linear coupled PDE system for which the density and the velocity are the unknowns. We use the finite element method to discretise the domain spatially. For stability reasons, we use Taylor-Hood elements \cite{girault_1979}. Moreover, it is recommended in \cite{rannacher_2000} to use quadrangle elements in 2D (respectively hexahedron in 3D) rather than triangle ones (respectively tetrahedron). Thus, we will use Qua4/Qua9 elements for the 2D flow problems and Hex8/Hex27 elements for the 3D problems.

The discretised unknowns $\pV$ and $\vV$ are gathered in the vector of unknows $\textbf{X} = \left(\pV \quad \vV \right)\Transpose$. The time derivatives for the unknowns are written as $\vVdot$ and $\pVdot$, and the right-hand side containing the external loadings is written $\textbf{F}$. Therefore, the matrix form of the problem is \eqref{eq:matrix_system}:
\begin{equation}\label{eq:matrix_system}
 \matD{M}(\vV,\pV) \mathbf{\dot{X}} + \matD{K}(\vV,\pV) \mathbf{X} = \mathbf{F}
\end{equation}

Equation \eqref{eq:matrix_system} is in a semi-discretised form, as the finite element method focuses only on the space domain $\dom$. The time domain $I = [0, T]$ is discretised in $N_t$ time steps, and we use a backward Euler scheme for the temporal integration.

\section{The LATIN-PGD solver for compressible flows}\label{sec:latin-pgd}

This section describes the LATIN-PGD algorithm for Newtonian compressible flows. The LATIN-PGD solver is an iterative non-incremental algorithm. It means that at each iteration of the solver, one gets a solution defined on the whole space-time domain. Two variants exist for this: the functional formulation and the internal variables formulation. For the sake of saving in computation times and memory, it is recommended in \cite{scanff_2021} to employ the functional formulation. In \autoref{sec:ref_pb}, we pointed out that flow models rely on conservation principles. Moreover, as underlined in the constitutive relations, the problem is coupled. The unknowns are the velocity and the density, which can be directly related to the pressure. Therefore, it is welcome to take inspiration from the multiphysics solver presented in \cite{dureisseix_2003} and adapted in a functional formulation in \cite{wurtzer_2024}.

\subsection{Reformulation of the problem}

The LATIN-PGD method relies on the separation of difficulties, which involves dividing all equations into two groups called $\Ad$ and $\Gam$. The first set of equations comprised in $\Ad$ should be linear and decoupled, and possibly global. The second group $\Gam$ is composed of all local equations, possibly coupled and non-linear. The reference solution is thus at the intersection of both sets of equations. We now clarify the equations that compose both sets $\Ad$ and $\Gam$. We introduce new variables to satisfy the properties of both groups stated above.

For the density part, we introduce the variables $\Z$, $\W$ and $\q$ and reformulate the continuity equation \eqref{eq:continuity}. The resulting equations \eqref{eq:global_eq_density} are in $\Ad$:
\begin{equation}\label{eq:global_eq_density}
\left\{
\begin{aligned}
    & \Z = \grad{\p} \quad \text{over}\: \Ixdom \\
    &  -\divergVec{\W} = \q \quad\text{over}\: \Ixdom
\end{aligned}
\right.
\end{equation}
For the velocity part, we introduce the variables $\Eps$ and $\Gamacc$ and reformulate the conservation of momentum \eqref{eq:momentum}. The resulting equations \eqref{eq:global_eq_velocity} in $\Ad$:
\begin{equation}\label{eq:global_eq_velocity}
\left\{
\begin{aligned}
    & \Eps = \frac{1}{2}(\gradVec{\V} + \gradVec{\V} \Transpose{})\quad\text{over}\: \Ixdom \\
    & \divergTens{\Sig} + \fd = \Gamacc \quad \text{over}\: \Ixdom
\end{aligned}
\right.
\end{equation}
We also include in $\Ad$ all the boundary conditions \eqref{eq:BCs_LATIN}:
\begin{equation}\label{eq:BCs_LATIN}
\left\{
\begin{aligned}
    & \p = \pd \quad\text{over} \: \bord{\rho}\times\I\\
    & \V = \Vd \quad\text{over} \: \bord{v}\times\I \\
    & \Sig \normal  = \Fd \quad \text{over}\: \bord{F}\times\I
\end{aligned}
\right.
\end{equation}

$\Gam$ gathers all equations that can be interpreted as constitutive ones \eqref{eq:gam_constitutive_rel}:
\begin{equation}\label{eq:gam_constitutive_rel}
\left\{
\begin{array}{llll}
    \Sig &= \Viscous \mathbin{:} \Eps - f(\p) \Id \\
    \Gamacc &= \p\derivFlu{\V}{t} = \p\displaystyle\derivPart{\V}{t} + \p\V \gradVec{(\V)} \\
    \W &= \p\V \\
    \q &= \displaystyle\derivPart{\p}{t}
\end{array}
\right.
\end{equation}
Those constitutive relations include time derivatives. Therefore, initial conditions \eqref{eq:IC_density} and \eqref{eq:IC_velocity} are also included in $\Gam$.

\subsection{Iterative algorithm with two alternated search directions}

As explained in the previous subsection, the reference solution of the flow problem satisfies:
\begin{equation}
 \sref = \Ad \cap \Gam
\end{equation}

The algorithm iterates between $\Ad$ and $\Gam$ until finding a good enough approximation of $\sref$. For that, each iteration of the algorithm is composed of two stages:
\begin{itemize}
    \item A \textbf{local stage}: from a known solution $\linn{\s} \in \Ad$, we compute a solution $\loc{\s} \in \Gam$ using a search direction $\Hup$;
    \item A \textbf{global stage}: from a known solution $\loc{\s} \in \Gam$, we compute a solution $\lin{\s} \in \Ad$ using a search direction $\Hdo$.
\end{itemize}

The search directions $\Hup$ and $\Hdo$ linking quantities at the local and global stages are given by:

\begin{equation}\label{eq:Hup}
\Hup \equiv \left\{
\Delta \s = (\Delta \sV, \Delta \sP) \quad \left|
\begin{array}{lllll}
    \Delta \Sig + \Hes^+ \mathbin{:} \Delta\Eps &= \tens{0} \\
    \Delta \Gamacc + \Hvg^+ \Delta\V &= \vect{0} \\
    \Delta \W + \Hzw^+ \Delta\Z &= \vect{0} \\
    \Delta \q + \Hpq^+ \Delta\p &= 0 \\
\end{array}\right.
\right\}
\end{equation}
\begin{equation}\label{eq:Hdo}
\Hdo \equiv \left\{\Delta \s = (\Delta \sV, \Delta \sP) \quad \left|
\begin{array}{lllll}
    \Delta \Sig - \Hes^- \mathbin{:} \Delta\Eps &= \tens{0} \\
    \Delta \Gamacc - \Hvg^- \Delta\V &= \vect{0} \\
    \Delta \W - \Hzw^- \Delta\Z &= \vect{0} \\
    \Delta \q - \Hpq^- \Delta\p &= 0 \\
\end{array}\right.
\right\}
\end{equation}
where  $\Delta\square$ is the difference of the quantity $\square$ between two consecutive stages. The choice of the search directions greatly influences the convergence rate of the algorithm \cite{scanff_2022}. As suggested in \cite{neron_2004,foulatier_2025}, it is relevant to choose here $\Hup$ and $\Hdo$  based on the constitutive laws \eqref{eq:gam_constitutive_rel} to optimise the convergence rate. Thus:
\begin{align*}
    &\Hes =\Hes^- =\Hes^+ = \Viscous \\
    &\Hvg = \Hvg^- = \Hvg^+ = \frac{1}{t_v}\Id\\
    &\Hzw =\Hzw^- =\Hzw^+ = -\frac{L_c^2}{T}\Id \\
    &\Hpq =\Hpq^- =\Hpq^+ = \frac{1}{t_\p}
\end{align*}

$L_c$ designates a characteristic length of the geometry, and $T$ corresponds to the final time step of the simulation. $t_v$ and $t_\p$ are respectively characteristic times for the velocity and the density parts that can be estimated with the demonstration in \cite{foulatier_2025}.

\subsection{Practical developments}

In this section, we provide a more detailed explanation of the operations required at each step of the LATIN-PGD algorithm. Before performing its iterations, the algorithm begins with an initialisation. It consists of computing two fields $\pVO$ and $\vVO$, both belonging to $\Ad$. Details concerning the problem solved at initialisation will be given in \autoref{sec:global_stage}.

\subsubsection{Local stage}

Knowing the solution from the previous iteration $\linn{\s} \in \Ad$, one computes the solution $\loc{\s} \in \Gam$. The search direction $\Hup$ yields:
\begin{subequations}
\begin{empheq}[left=\empheqlbrace]{align}
    &\loc{\Sig} + \Hes \mathbin{:}\loc{\Eps} = \underbrace{\linn{\Sig} + \Hes\mathbin{:}\linn{\Eps}}_{\displaystyle\linn{\Atens}} \\
    &\loc{\Gamacc} + \Hvg \loc{\V} =  \underbrace{\linn{\Gamacc} + \Hvg\linn{\V}}_{\displaystyle\linn{\betavec}} \\
    &\loc{\W} + \Hzw \loc{\Z} =  \underbrace{\linn{\W} + \Hzw\linn{\Z}}_{\displaystyle\linn{\deltavec}} \\
    &\loc{\q} + \Hpq \loc{\p} = \underbrace{\linn{\q} + \Hpq\linn{\p}}_{\displaystyle\linn{\gammascal}}
\end{empheq}
\end{subequations}

Using the constitutive equations, one gets:
\begin{subequations}
\begin{empheq}[left=\empheqlbrace]{align}
    \label{eq:coupled1} &(\Viscous + \Hes)\mathbin{:} \loc{\Eps} - f\left(\loc{\p}\right)\Id =  \linn{\Atens} \\
    \label{eq:coupled2} &\loc{\p}\derivPart{\loc{\V}}{t} + \loc{\p}\loc{\V} \gradVec{\loc{\V}} + \Hvg\loc{\V} = \linn{\betavec} \\
    \label{eq:coupled3} &\loc{\p}\loc{\V} + \Hzw \loc{\Z} = \linn{\deltavec} \\
    \label{eq:coupled4} &\derivPart{\loc{\p}}{t} + \Hpq\loc{\p} = \linn{\gammascal}
\end{empheq}
\end{subequations}

\eqref{eq:coupled4} is an uncoupled ordinary differential equation and is solved using a backward Euler scheme. Then, injecting $\loc{\p}$ in \eqref{eq:coupled1}, one directly gets $\loc{\Eps}$. Knowing $\loc{\p}$, the only unknown in \eqref{eq:coupled2} is $\loc{\V}$. However, \eqref{eq:coupled2} is a non-linear equation. To solve it, we can use a Newton-Raphson algorithm. Finally, assuming $\loc{\p}$ and $\loc{\V}$ are known, the density gradient $\loc{\Z}$ is directly computed with \eqref{eq:coupled3}. As all of these equations are local, they can be solved independently for each integration point.

The study focuses here on laminar flow, which means that the non-linear term in \eqref{eq:coupled2} can be neglected. Therefore, the resolution of the equation just requires classical methods to solve an ordinary differential equation of order 1. The system solved at the coupled stage is simplified to:
\begin{subequations}
\begin{empheq}[left=\empheqlbrace]{align}
    \label{eq:coupled1_bis} &(\Viscous + \Hes)\mathbin{:} \loc{\Eps} - f\left(\loc{\p}\right)\Id =  \linn{\Atens} \\
    \label{eq:coupled2_bis} &\loc{\p}\derivPart{\loc{\V}}{t} + \Hvg\loc{\V} = \linn{\betavec} \\
    \label{eq:coupled3_bis} &\loc{\p}\loc{\V} + \Hzw \loc{\Z} = \linn{\deltavec} \\
    \label{eq:coupled4_bis} & \derivPart{\loc{\p}}{t} + \Hpq\loc{\p} = \linn{\gammascal}
\end{empheq}
\end{subequations}

Finally, one deduces the dual quantities with the search direction $\Hup$ as:
\begin{subequations}
\begin{empheq}[left=\empheqlbrace]{align}
        &\loc{\Sig} = \linn{\Atens} - \Hes\mathbin{:}\loc{\Eps}\\
        &\loc{\Gamacc} = \linn{\betavec} - \Hvg\loc{\V}\\
        &\loc{\W} = \linn{\deltavec} - \Hzw\loc{\Z}\\
        &\loc{\q} = \linn{\gammascal} - \Hpq\loc{\p}
\end{empheq}
\end{subequations}

\subsubsection{Global stage}\label{sec:global_stage}

The search direction $\Hdo$ yields:
\begin{subequations}
\begin{empheq}[left=\empheqlbrace]{align}
    \label{eq:decoupled1} & \lin{\Sig} = (\underbrace{\loc{\Sig} - \Hes\mathbin{:}\loc{\Eps}}_{\displaystyle\loc{\Atens}}) + \Hes\mathbin{:}\lin{\Eps}\\
    \label{eq:decoupled2} & \lin{\Gamacc} = \underbrace{(\loc{\Gamacc} - \Hvg\loc{\V})}_{\displaystyle\loc{\betavec}} + \Hvg\lin{\V} \\
    \label{eq:decoupled3} & \lin{\W} = \underbrace{(\loc{\W} - \Hzw\loc{\Z})}_{\displaystyle\loc{\deltavec}} + \Hzw\lin{\Z} \\
    \label{eq:decoupled4} & \lin{\q} = \underbrace{(\loc{\q} - \Hpq\loc{\p})}_{\displaystyle\loc{\gammascal}} + \Hpq\lin{\p}
\end{empheq}  
\end{subequations}

Let us first focus on the weak formulation of the continuity equation in which equations \eqref{eq:decoupled3} and \eqref{eq:decoupled4} are injected. The global problem to solve is then:

Find $\lin{\p} \in \spaceP$ such that $\forall t \in \I$: 
\begin{equation}\label{eq:weak_form_rho}
\begin{split}
    -\intO{\grad{\lin{\p}}\cdot\Hzw\grad{\Ptest}}+\intO{\lin{\p}\Hpq\Ptest}=&\intO{\loc{\deltavec}\cdot\grad{\Ptest}} \\
        &-\intO{\loc{\gammascal}\Ptest} \qquad \forall \Ptest \in \spacePo
\end{split}
\end{equation}
At each iteration, we compute a correction $\Delta\lin{\p}$ such that: $\lin{\p} = \linn{\p}+\Delta\lin{\p}$. Therefore, the problem transforms to:
\begin{equation}\label{eq:weak_form_rho_bis}
\begin{split}
   -\intO{\grad{(\Delta\lin{\p})}\cdot\Hzw\grad{\Ptest}} &+\intO{\Delta\lin{\p}\Hpq\Ptest}= \intO{\loc{\deltavec}\cdot\grad{\Ptest}}-\intO{\loc{\gammascal}\Ptest}\\
    & \intO{\grad{\linn{\p}}\cdot\Hzw\grad{\Ptest}}-\intO{\linn{\p}\Hpq\Ptest} \qquad \forall \Ptest \in \spacePo
\end{split}
\end{equation}

At the initialisation, one takes $\hat{\deltavec}_0 = \vect{0}$ and $\hat{\gammascal}_0 = 0$. We seek an initial guess of $\pV$ by computing the finite element problem \eqref{eq:p_FE_init} corresponding to \eqref{eq:weak_form_rho_bis}:
\begin{equation}\label{eq:p_FE_init}
\left\{
    \begin{aligned}
        \Mpp{H} \pVO(t) &= \fp{f}(t) \quad \text{over}\:\dom\\
        \pVO(t) &= \pVd(t) \quad \text{over}\: \bord{\p}
    \end{aligned}
\right.
\end{equation}
where $\fp{f}$ represents the external loadings. The matrix $\Mpp{H}$ is bound to $\Hzw$ and $\Hpq$ in such way:
\begin{equation}
    \Mpp{H} = -\Hzw \Mpp{K} + \Hpq \Mpp{C}    
\end{equation}
where $\Mpp{C}$ and $\Mpp{K}$ are the finite element matrices respectively corresponding to the integral over the domain of shape functions products and the integral of their 1st-order derivatives products.

At any iteration $n>0$, the discrete form for the weak formulation is given by \eqref{eq:decoupled_density}:
\begin{equation}\label{eq:decoupled_density}
\left\{
    \begin{aligned}
        \Mpp{H} \linCorP(t) &= \fp{g}(t) \quad \text{in}\:\dom\\
        \linCorP(t) &= \vecD{0} \quad \text{on}\: \bord{D}
    \end{aligned}
\right.
\end{equation}
with $\fp{g} = \Fdc + \Fgc - \Mpp{H}(\linShort{\pV}-\pVO)$ where $\Fdc$ and $\Fgc$ refer to the right-hand side computed with the local stage quantities.

The same procedure is applied to obtain the initialisation $\vVO$ and correction $\linCorV$ of the velocity. Using equations \eqref{eq:decoupled1} and \eqref{eq:decoupled2} yields the following weak form:

Find $\lin{\V} \in \spaceV$ such that $\forall t \in \I$: 
\begin{equation}\label{eq:weak_form_v}
\begin{split}
    \intO{\Eps(\Delta\lin{\V})\mathbin{:}\Hes\mathbin{:}\Eps(\Vtest)} + & \intO{\scalprod{\Delta\lin{\V}}{\Hvg\Vtest}}=\\
    & - \intO{\loc{\Atens}:\Eps(\Vtest)}- \intO{\scalprod{\loc{\betavec}}{\Vtest}}\\
    &+\intO{\scalprod{\vect{b}}{\Vtest}} + \intS{N}{\scalprod{\Fd}{\Vtest}}\\
    & -\intO{\Eps(\linn{\V})\mathbin{:}\Hes\mathbin{:}\Eps(\Vtest)} - \intO{\scalprod{\linn{\V}}{\Hvg\Vtest}} \qquad \forall \Vtest \in \spaceVo
\end{split}
\end{equation}

At the initialisation, one takes $\hat{\Atens}_0 = \tens{0}$ and $\hat{\betavec}_0 = \vect{0}$ and we compute the initial guess $\vVO$ by solving:
\begin{equation}
\left\{
    \begin{aligned}
        \Mvv{H} \vVO(t) &= \fv{f}(t) \quad \text{over}\:\dom\\
        \vVO(t) &= \vVd(t) \quad \text{over}\: \bord{v}
    \end{aligned}
\right.
\end{equation}
where $\fv{f}$ represent the external loadings. The matrix $\Mvv{H}$ is bound to $\Hes$ and $\Hvg$ in such way:
\begin{equation}
    \Mvv{H} = \Hes \Mvv{K} + \Hvg \Mvv{C}    
\end{equation}
$\Mvv{K}$ and $\Mvv{C}$ are defined in the same way as $\Mpp{K}$ and $\Mpp{C}$ but correspond to the velocity elements.

At any iteration $n>0$, the discrete form for the weak formulation is given by \eqref{eq:decoupled_velocity}:
\begin{equation}\label{eq:decoupled_velocity}
\left\{
    \begin{aligned}
        \Mvv{H} \linCorV(t) &= \fv{g}(t) \quad \text{in}\:\dom\\
        \linCorV(t) &= \vecD{0} \quad \text{on}\: \bord{v}
    \end{aligned}
\right.
\end{equation}
with $\fv{g} = \FAc + \Fbc - \Mvv{H}(\linShort{\vV}-\vVO)$ where $\FAc$ and $\Fbc$ refer to the right-hand side computed with the local stage quantities.

\subsubsection{Proper Generalised Decomposition within the LATIN framework}

The strength of the method lies in solving the decoupled stage using \textit{Proper Generalised Decomposition} (PGD). Thus, it allows us to express the primal fields as a combination of time and space functions. After $n$ iterations of the algorithm, the density $\pV$ and the velocity $\vV$ are expressed as follows:
\begin{equation}\label{eq:primal_fields_PGD}
    \left\{
    \begin{aligned}
    \linn{\vV} (t) &= \vVO(t) + \sumPGD{i}{k} \TimeModeV{i}{}(t) \SpaceModeV{i} \quad k\leqslant n\\
    \linn{\pV} (t) &= \pVO(t) + \sumPGD{i}{l} \TimeModeP{i}(t) \SpaceModeP{i} \quad l\leqslant n
    \end{aligned}
    \right.
\end{equation} 

At each iteration and for each primal variable, one seeks a new pair of modes in the form given in \eqref{eq:gen_modes_PGD} if needed:
\begin{equation}\label{eq:gen_modes_PGD}
    \left\{
    \begin{aligned}
    \linCorV (t) &= \TimeModeV{}{}(t) \SpaceModeV{}\\
    \linCorP (t) &= \TimeModeP{}(t) \SpaceModeP{}
    \end{aligned}
    \right.
\end{equation}

The decoupled stage always begins with an update step that changes all time functions at a significantly lower cost compared to generating a new spatial mode. If the update step has sufficiently improved the solution, we do not perform any mode generation afterward. For more practical details on the use of PGD at the LATIN decoupled stage, one can refer to \cite{heyberger_2012}.

\subsection{Convergence indicator}

It is necessary to define a consistent criterion to stop the iterations. We define an error indicator quantifying the gap between two successive solutions $\loc{\s}$ and $\lin{\s}$ at the end of an iteration. When we reach a sufficiently small difference, the stopping criterion claims that the approximation of the reference solution is good enough. 

We introduce two error indicators, one related to density and the other to velocity. For each of these indicators, we introduce a specific energetic norm.

The density error indicator is given by:
\begin{equation}\label{eq:ind_pressure}
 {\etaP} ^2 = \frac{\norme{\lin{\sP}-\loc{\s}}{\Hpq}^2}{\frac{1}{2}(\norme{\lin{\sP}}{\Hpq}^2 + \norme{\loc{\s}}{\Hpq}^2)} \quad \text{with} \norme{\sP}{\Hpq}^2 = \intTO{\p \Hpq \p}
\end{equation}

The velocity error indicator is given by:
\begin{equation}\label{eq:ind_vel}
 {\etaV} ^2 = \frac{\norme{\lin{\sV}-\loc{\s}}{\Hes}^2}{\frac{1}{2}(\norme{\lin{\sV}}{\Hes}^2 + \norme{\loc{\s}}{\Hes}^2)} \quad \text{with} \norme{\sV}{\Hes}^2 = \intTO{\Eps : \Hes : \Eps}
\end{equation}
One provides a tolerance threshold denoted $\eta_c$ and the algorithm ends when $\eta = \text{max}\left(\etaV, \etaP \right)~<~\eta_c$. Therefore, one can stop computing the global stage of one part if it has converged faster than the other.

\subsection{Summary of the algorithm}

\begin{figure}[h!]
    \centering
    \includegraphics[width=15cm]{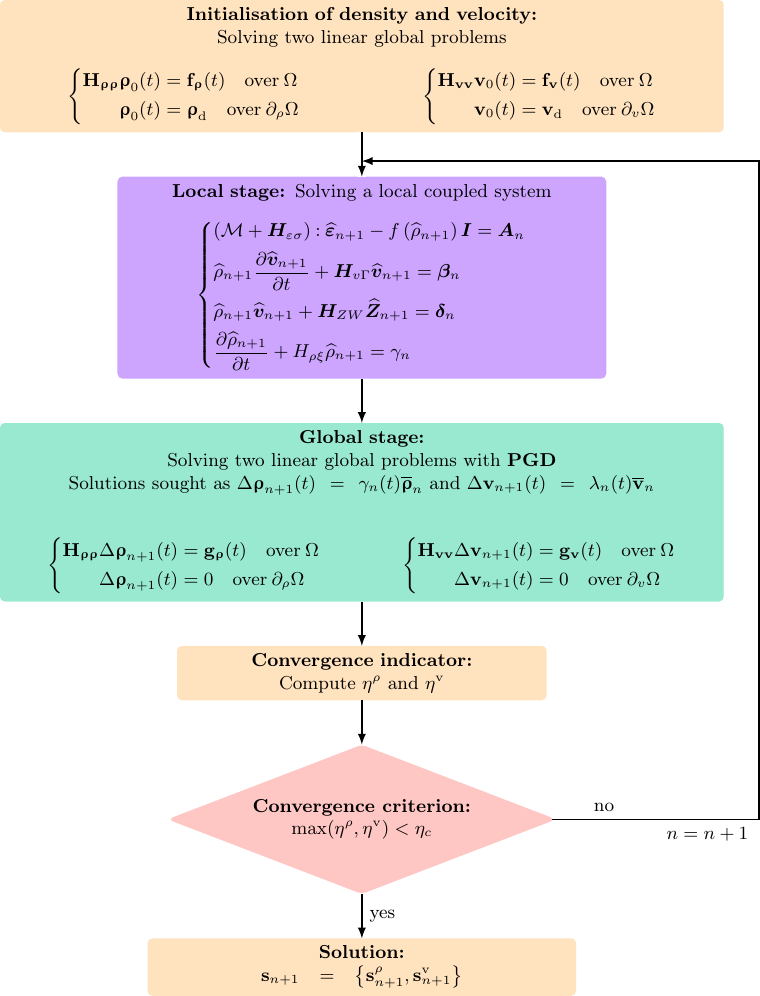}
    \caption{Flowchart of the LATIN-PGD algorithm for compressible flows}
    \label{fig:algo_LATIN}
\end{figure}

The flowchart in \autoref{fig:algo_LATIN} sums up the LATIN-PGD algorithm for Newtonian compressible flows.

\section{Results}\label{sec:results}

\subsection{Laminar flow in a channel}\label{sec:channel}

We first apply the method to a simple example involving a channel subjected to a pressure difference between its inflow and outflow boundaries. We will first compare the results to a given analytical solution at steady state. Then, we will display the results for the entire space-time domain and provide further details about the PGD modes obtained during the simulation.

\subsubsection{Presentation of the problem}\label{sec:channel_pres}

\begin{figure}[ht!]
    \centering
    \includegraphics[width=15cm]{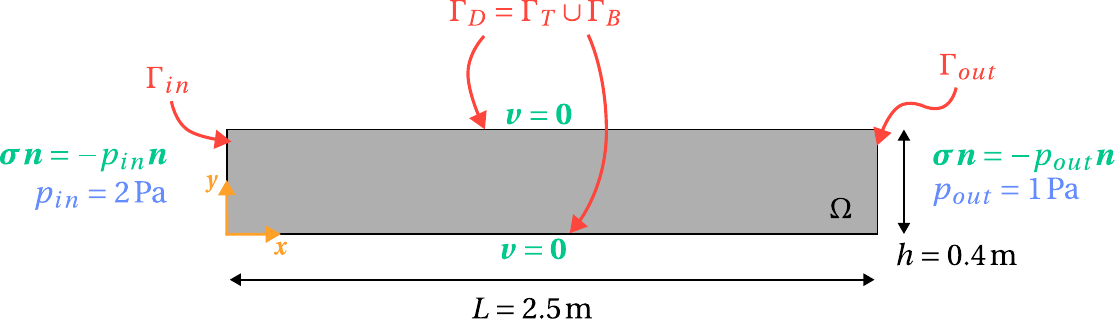}
    \caption{Scheme of the studied problem}
    \label{fig:scheme_channel}
\end{figure}

The geometry and boundary conditions of the problem are represented in \autoref{fig:scheme_channel}, where the $x$-axis is oriented along the channel axis. We consider a pressure $p_{in} = \qty{2}{\pascal}$ on the inflow boundary $\Gamma_{in}$ and a pressure $p_{out} = \qty{1}{\pascal}$ on the outflow boundary $\Gamma_{out}$. On the top and bottom boundaries $\Gamma_D$, we consider a homogeneous Dirichlet condition for the velocity part (no-slip boundary condition). Moreover, for the velocity part, one must consider Neumann boundary conditions due to the inflow and outflow pressures, respectively, on $\Gamma_{in}$ and $\Gamma_{out}$.

\autoref{tab:para_mat_channel} summarises the material properties. One notices a considerable difference between the numerical values of the first and second viscosities, as shown in \cite{liebermann_1949,rosenhead_1954}.

\begin{table}[ht!]
    \centering
    \caption{Fluid material parameters used for the laminar flow in a channel problem}
    \label{tab:para_mat_channel}
    \begin{tabular}{c c}
    \toprule
    Parameters & Value (SI units) \\
    \hline
        Dynamic viscosity $\mu$ & $\qty{1}{\kilo\gram\per\meter\per\second}$ \\
        Second viscosity $\lambda$ & $\qty{1e3}{\kilo\gram\per\meter\per\second}$ \\
        Reference temperature $T_0$ & $\qty{293}{\kelvin}$\\
        Universal gas constant $R$ & $\qty{8.314}{\joule\per\kelvin\per\mol}$\\
        Molar mass $M$ & $\qty{28.9645e-3}{\kilogram\per\mol}$\\
    \bottomrule
    \end{tabular}
\end{table}

For the initial conditions, one takes:
\begin{subequations}
\begin{empheq}[left=\empheqlbrace]{align}
    \label{eq:velocity_init}& \V(t=0, \M) = \vect{0} & \forall \M \in \dom \\
    \label{eq:density_init}& \p(t=0, \M) = \po = \frac{p_0}{rT_0} & \forall \M \in \dom \\
    \label{eq:pressure_init}& p(t=0, \M) = p_0 = \qty{1}{\pascal} & \forall \M \in \dom
\end{empheq}
\end{subequations}

We use a structured mesh with 128 elements along the length of the channel and 16 elements along the height of the channel, totaling 2048 Qua4/Qua9 elements. It results in 16,962 DoFs for the velocity and 2,193 DoFs for the density.

\subsubsection{Results at steady state}

We first focus on the results at steady state for which an analytical solution is known. Under the assumption of stationarity, one has to solve the Stokes equations and analytically finds:
\begin{subequations}
\begin{empheq}[left=\empheqlbrace]{align}
    \label{eq:velocity_ana}& \V(\M) = \frac{\left((h/2)^2 - y^2\right)}{2\mu}\frac{\lvert \Delta p \rvert}{L} \vect{e}_x & \forall \M \in \dom\\
    \label{eq:pressure_ana}& p(\M) = p_{in} + \frac{\Delta p}{L}x & \forall \M \in \dom
\end{empheq}
\end{subequations}
with $\Delta p = p_{out} - p_{in}$. The velocity is then maximal on the mean line and equal to \eqref{eq:vmax_ana}:
    \begin{equation}\label{eq:vmax_ana}
 v_{max} = \cfrac{(h/2)^2}{2\mu}\cfrac{\lvert \Delta p \rvert}{L}
    \end{equation}

For the studied problem, we have the numerical value $v_{max} = \qty{8e-3}{\meter\per\second}$.

\begin{figure}[ht!]
    \centering
    \begin{subfigure}[b]{0.9\textwidth}
         \centering
         \includegraphics[width=10.5cm]{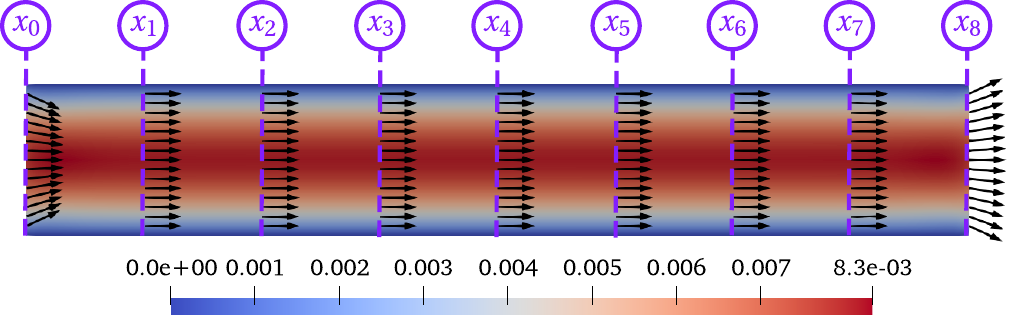}
         \caption{Velocity [m/s]: $x$-axis component}
         \label{fig:vel_stationary}
     \end{subfigure}
     \hfill
     \begin{subfigure}[b]{0.9\textwidth}
         \centering
         \includegraphics[width=10cm]{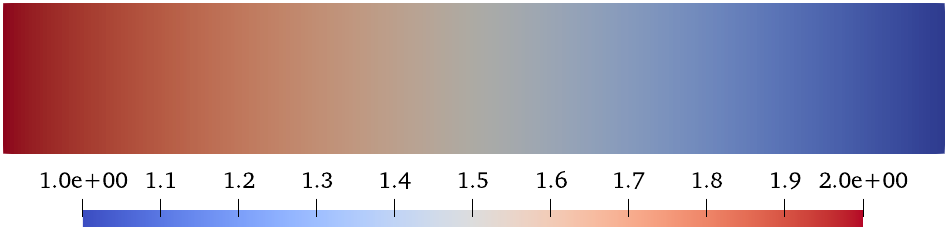}
         \caption{Pressure [Pa]}
         \label{fig:pressure_stationary}
     \end{subfigure}
    \caption{Fields at steady state}
    \label{fig:fields_stationary}
\end{figure}

The velocity and pressure fields obtained with the LATIN-PGD method are depicted in \autoref{fig:fields_stationary}.

We qualitatively recover the abovementioned properties for both fields. The pressure field is uniformly decreasing along the $x$-axis, and the velocity field is constant in each channel cross-section, reaching its maximal value on the mean line. On \autoref{fig:vel_stationary}, one sees that the $y$ component of the velocity field is equal to zero almost everywhere except at the boundaries. We indeed expect the flow to be oriented towards the mean line at the inflow and outwards at the outflow. 

\begin{figure}[ht!]
    \centering
    \begin{subfigure}[b]{0.99\textwidth}
         \centering
         \includegraphics[width=10cm]{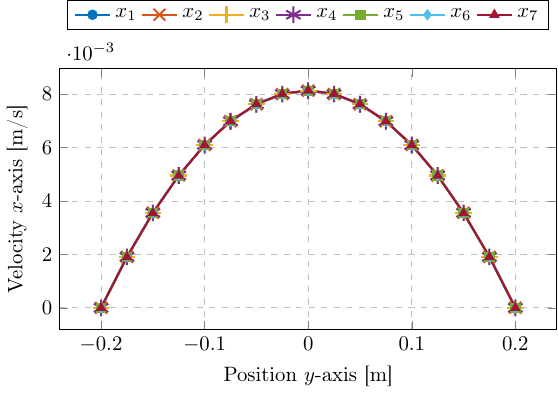}
         \caption{Velocity field in different sections of the channel}
         \label{fig:velocity_7_sections}
     \end{subfigure}
     \hfill
     \begin{subfigure}[b]{0.99\textwidth}
         \centering
         \includegraphics[width=10cm]{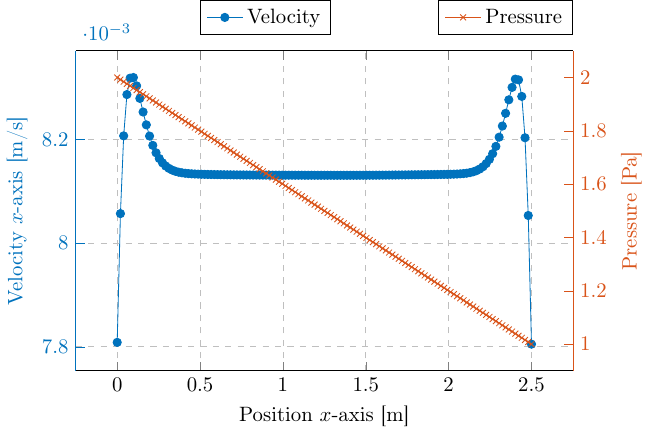}
         \caption{Velocity and pressure along mean line}
         \label{fig:fields_avg_line}
     \end{subfigure}
    \caption{Velocity and pressure results at steady state}
    \label{fig:plots_steady_state}
\end{figure}

On \autoref{fig:velocity_7_sections}, we plot the value of velocity along the $y$-axis at different sections $x_k$ such that $x_k~=~k~\times~\qty{0.3125}{\meter}, \: k \in \mathbb{N}$. First, one can check that the velocity is almost constant in each section of the channel. Secondly, the velocity field has exactly the shape of a parabolic function. Finally, the maximal value is very close to $\qty{8e-3}{\meter\per\second}$, which is the value expected. On \autoref{fig:fields_avg_line}, we depict the $x$ component of the velocity and the pressure field along the mean line of the channel. We clearly notice the linear decrease in pressure. The velocity is almost constant everywhere, except at $x=0$ and $x=L$, where boundary effects occur.

\begin{figure}[ht!]
    \centering
    \begin{subfigure}[b]{0.9\textwidth}
         \centering
         \includegraphics[width=10.5cm]{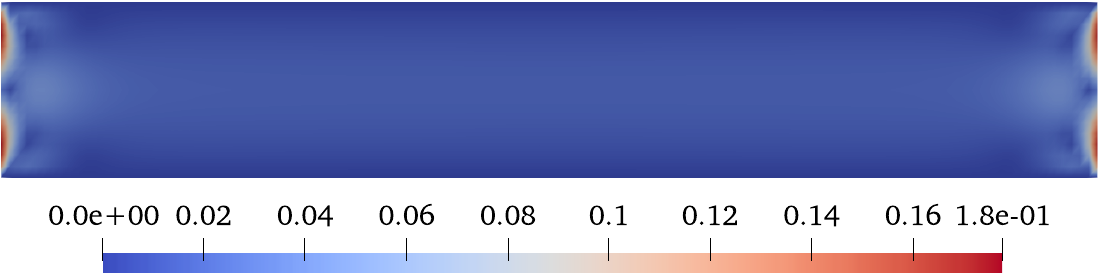}
         \caption{Relative error on velocity magnitude [-]}
         \label{fig:channel_err_vel}
     \end{subfigure}
     \hfill
     \begin{subfigure}[b]{0.9\textwidth}
         \centering
         \includegraphics[width=10cm]{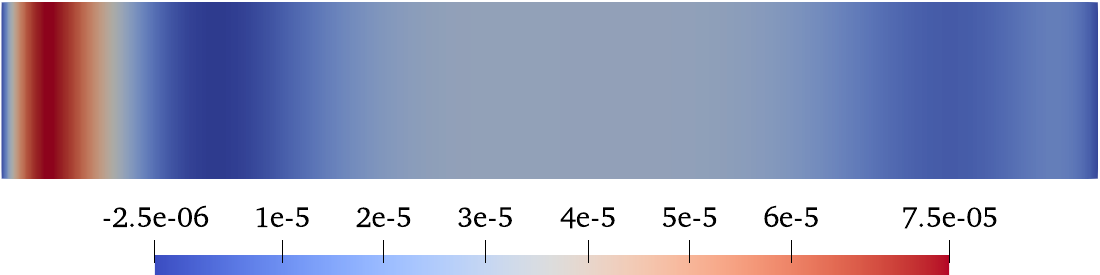}
         \caption{Relative error on pressure [-]}
         \label{fig:channel_err_pressure}
     \end{subfigure}
    \caption{Relative error between numerical and analytical solution}
    \label{fig:relative_errors}
\end{figure}

We plot on \autoref{fig:relative_errors} the relative error between the LATIN-PGD solution and the analytical one. 

We have a very low error rate, specifically in the pressure data. Close to the boundary, we remark slightly higher errors due to the error in the $y$ component of the velocity. We explain those errors by assuming an infinite channel in the analytical calculations.

\subsubsection{Results for the whole space-time domain}

The results presented in the previous section align with the solution of the Stokes equation when the steady state is reached. However, the Navier-Stokes equations include time-dependency. Thus, this section will focus on the results for the whole time domain $I = [0, 5] \: \qty{}{\milli\second}$. The time domain is divided into $100$ time steps, resulting in a time step size of $\Delta t = \qty{5e-2}{\milli\second}$.

\begin{figure}[ht!]
    \centering
    \begin{subfigure}[b]{0.49\textwidth}
         \centering
         \includegraphics[width=7cm]{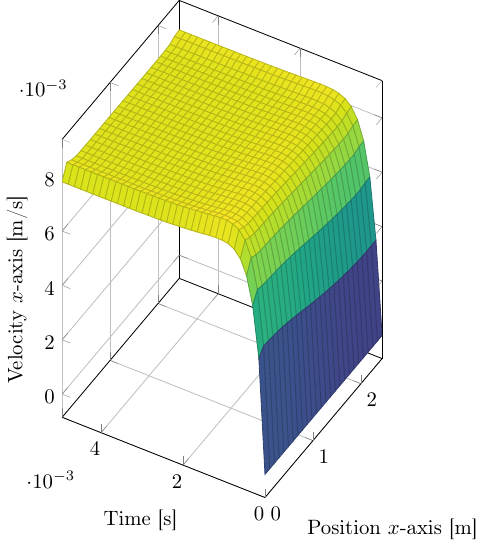}
         \caption{Velocity}
         \label{fig:surf_velocity}
     \end{subfigure}
     \hfill
     \begin{subfigure}[b]{0.49\textwidth}
         \centering
         \includegraphics[width=7cm]{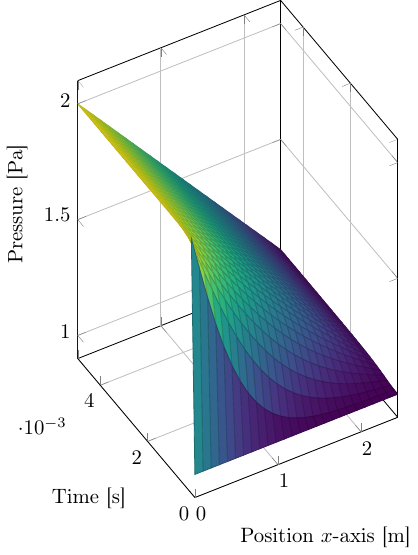}
         \caption{Pressure}
         \label{fig:surf_pressure}
     \end{subfigure}
    \caption{Velocity and pressure on the mean line along time}
    \label{fig:surface_plots}
\end{figure}

One can first examine the evolution of the velocity $x$-component and the pressure over time. We plot both quantities in \autoref{fig:surface_plots} as a function of position along the mean line and time.

\begin{figure}[ht!]
    \centering
    \includegraphics[width=8cm]{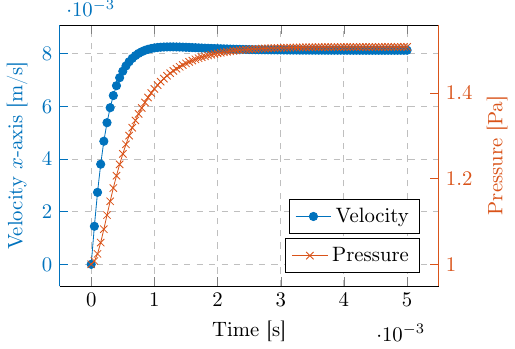}
    \caption{Evolution of pressure and velocity along time at point $(\qty{1.25}{},0)$}
    \label{fig:fields_middle_pt_vs_time}
\end{figure}

We observe that the solution now exhibits a time dependency and reaches its steady state within a few milliseconds. In order to have a closer look at the unsteady state, we display on \autoref{fig:fields_middle_pt_vs_time} velocity and pressure at the point $(\qty{1.25}{},0)$ (corresponding to the middle point of the channel). The unsteady state lasts from $\qty{0}{\milli\second}$ to approximatively $\qty{3}{\milli\second}$. The quantities at steady state fit then with the numerical values expected in the middle of the channel.

\begin{figure}[ht!]
    \centering
    \begin{subfigure}[b]{0.49\textwidth}
         \centering
         \includegraphics[width=7cm]{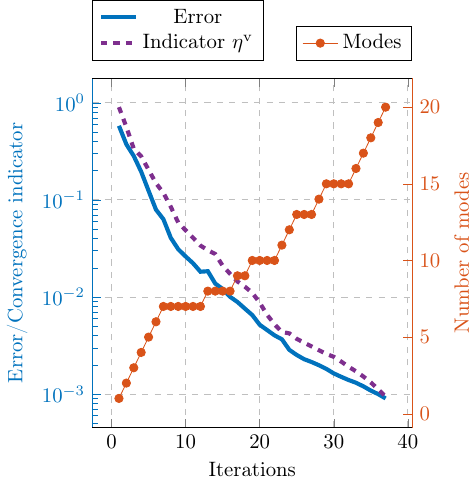}
         \caption{Velocity}
         \label{fig:latin_pgd_vel}
     \end{subfigure}
     \hfill
     \begin{subfigure}[b]{0.49\textwidth}
         \centering
         \includegraphics[width=7cm]{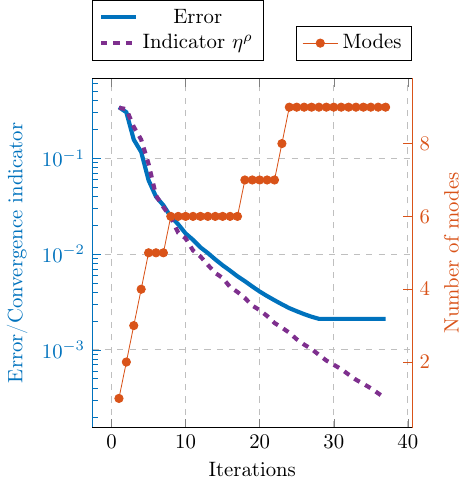}
         \caption{Density}
         \label{fig:latin_pgd_dens}
     \end{subfigure}
    \caption{Evolution of error, indicator, and number of modes along LATIN iterations for velocity and density parts}
    \label{fig:latin_pgd_channel}
\end{figure}

It is also important to discuss the modes generated by PGD. \autoref{fig:latin_pgd_channel} displays the convergence indicators $\etaV$ and $\etaP$ along iterations and the error towards a full-order solution obtained with a very low convergence threshold ($\eta_c = \qty{1e-8}{}$). Moreover, one can read the number of modes along the computation. We used the numerical values $t_v = T_f$ and $t_\rho = T_f/10$ to obtain following results.

\begin{figure}[ht!]
    \centering
    \begin{subfigure}[b]{0.49\textwidth}
         \centering
         \includegraphics[width=\textwidth]{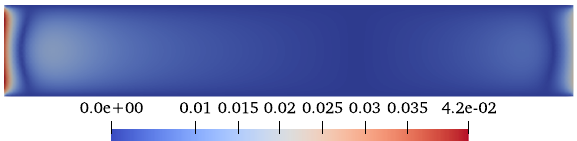}
         \caption{Mode 1 [m/s]}
         \label{fig:mode1}
     \end{subfigure}
     \hfill
     \begin{subfigure}[b]{0.49\textwidth}
         \centering
         \includegraphics[width=\textwidth]{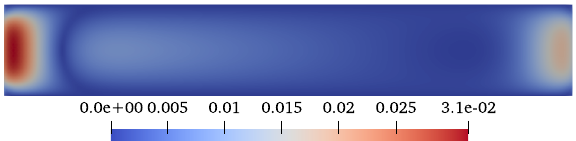}
         \caption{Mode 2 [m/s]}
         \label{fig:mode2}
     \end{subfigure}
     \hfill
     \begin{subfigure}[b]{0.49\textwidth}
         \centering
         \includegraphics[width=\textwidth]{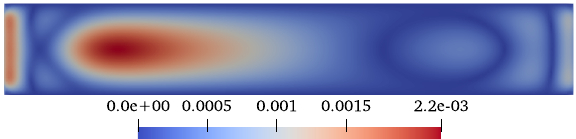}
         \caption{Mode 4 [m/s]}
         \label{fig:mode4}
     \end{subfigure}
     \hfill
     \begin{subfigure}[b]{0.49\textwidth}
         \centering
         \includegraphics[width=\textwidth]{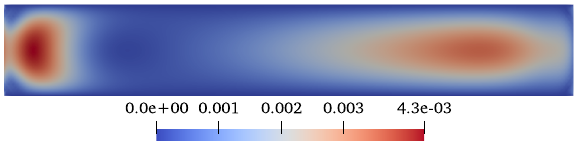}
         \caption{Mode 5 [m/s]}
         \label{fig:mode5}
     \end{subfigure}
    \caption{Some velocity modes at the last time step given by the LATIN-PGD algorithm}
    \label{fig:modes_pgd}
\end{figure}

The number of modes is higher for the velocity than for the density, which is expected, as the solution is more complex. To understand why we require more modes to describe the velocity, it is interesting to observe the first velocity modes generated (see \autoref{fig:modes_pgd}).

The first two modes mainly depict the velocity changes near the inflow and outflow boundaries. As we employ an on-the-fly reduction technique, it is logical that the first modes are not the optimal ones. Indeed, we initialise the algorithm with a solution far from the reference solution. The modes generated afterwards (modes 3 to 8) represent the evolution of velocity inside the channel. With these first eight modes, we qualitatively recover the shape of velocity from the reference solution. However, the error is still comprised between $5 \%$ and $10\%$. The following modes have a lower amplitude and locally improve the solution.

\subsection{Resolution of the 2D-1 benchmark}

\begin{figure}[ht!]
    \centering
    \includegraphics[width=12cm]{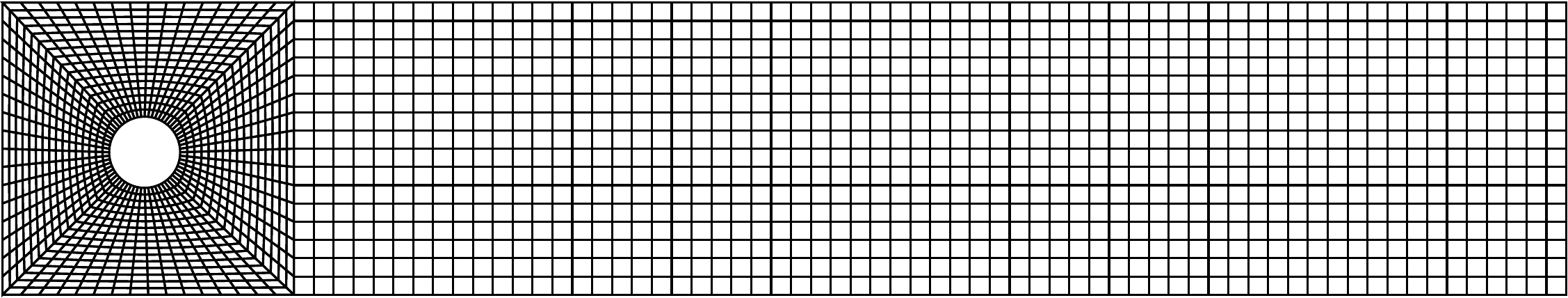}
    \caption{Mesh used for the 2D-1 benchmark}
    \label{fig:mesh_2D_benchmark}
\end{figure}

We now want to reproduce the 2D-1 benchmark for laminar flows. We consider a channel of length $L=\qty{2.2}{\meter}$ and height $h = \qty{0.41}{\meter}$, where $x$-axis is oriented along the channel main axis. As described in \cite{schafer_1996}, it includes a cylinder of radius $r=\qty{0.05}{m}$ whose origin is located at $(0.2, \: 0.2)\qty{}{\meter}$. The structured mesh used for this benchmark is given in \autoref{fig:mesh_2D_benchmark}. The mesh is divided into 2,368 quadrangle elements, resulting in 16,896 DoFs for the velocity and 2,176 DoFs for the density.

The material properties are the same as the ones given in \autoref{tab:para_mat_channel}.

\begin{figure}[ht!]
    \centering
     \begin{subfigure}[b]{0.99\textwidth}
         \centering
         \includegraphics[width=10cm]{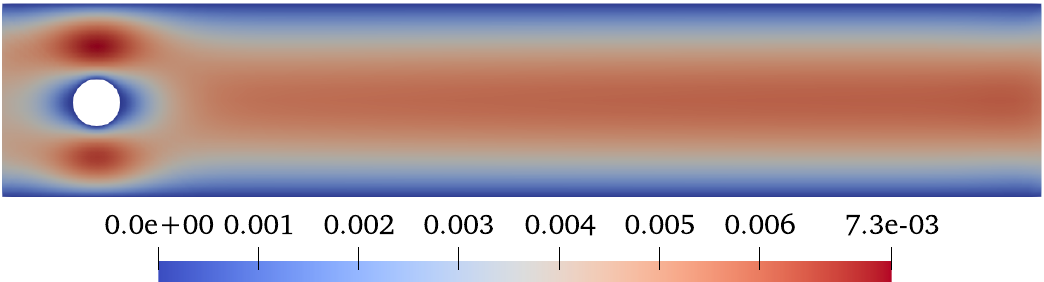}
         \caption{Velocity [m/s]}
         \label{fig:vel_benchmark}
     \end{subfigure}
     \hfill
     \begin{subfigure}[b]{0.99\textwidth}
         \centering
         \includegraphics[width=10cm]{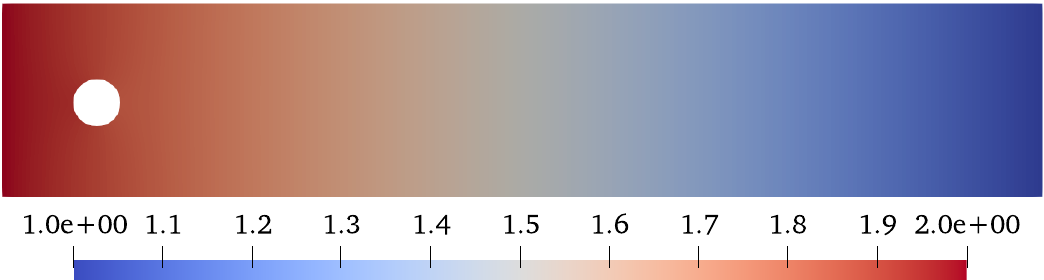}
         \caption{Pressure [Pa]}
         \label{fig:pressure_benchmark}
     \end{subfigure}
    \caption{Fields at last time step for the 2D-1 benchmark}
    \label{fig:fields_benchmark}
\end{figure}

\autoref{fig:fields_benchmark} displays the velocity and pressure fields at the last time step. We observe a maximum velocity at the cylinder boundaries, which is consistent with the literature \cite{schafer_1996}. Compared to other studies, we have perfect symmetry in the velocity fields between the left and right sides of the cylinder, since we vanish the nonlinear convective term entirely.

\begin{figure}[ht!]
    \centering
    \centering
    \includegraphics[width=8cm]{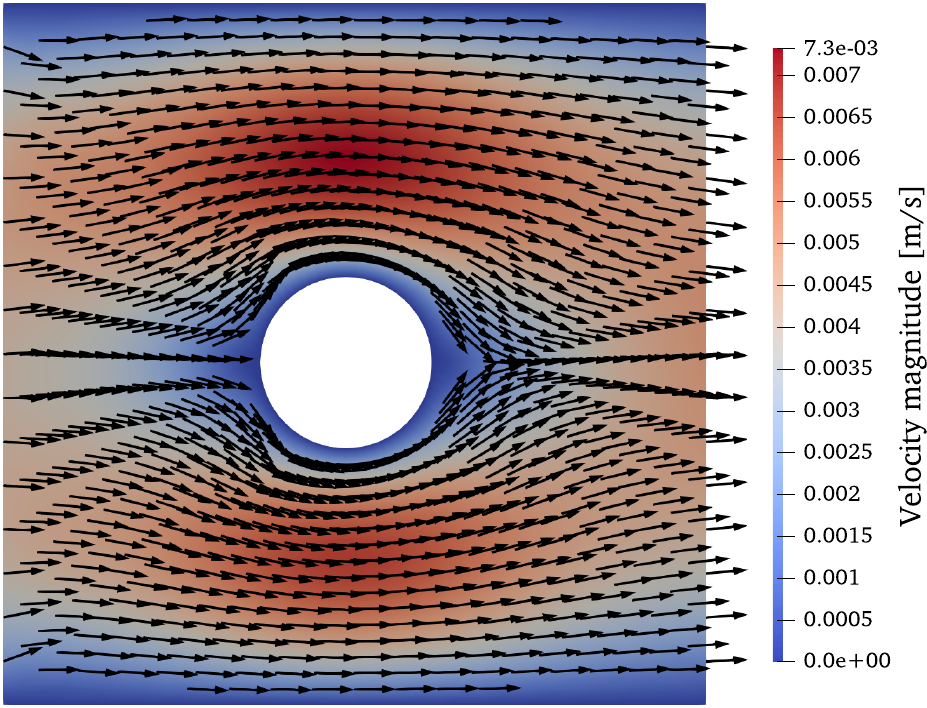}
    \caption{Velocity field and stream lines around the cylinder}
    \label{fig:zoom_velocity}
\end{figure}

The \autoref{fig:zoom_velocity} displays a zoom around the cylinder, enabling one to see the streamlines around the cylinder. The velocity streamlines are bypassing the cylinder in a regular way, which is relevant for laminar flows.

At the end of the computation, we have built a basis representing the velocity with 15 modes and another with 6 modes for the pressure. With such bases, we achieve a global relative error of $0.1\%$ for the velocity and $0.3\%$ for the pressure in 39 iterations. Despite the higher complexity of the benchmark, the convergence speed and basis size are very similar to those presented in \autoref{sec:channel}.

\subsection{Extension to 3D problems}

\begin{figure}[ht!]
    \centering
    \includegraphics[width=15cm]{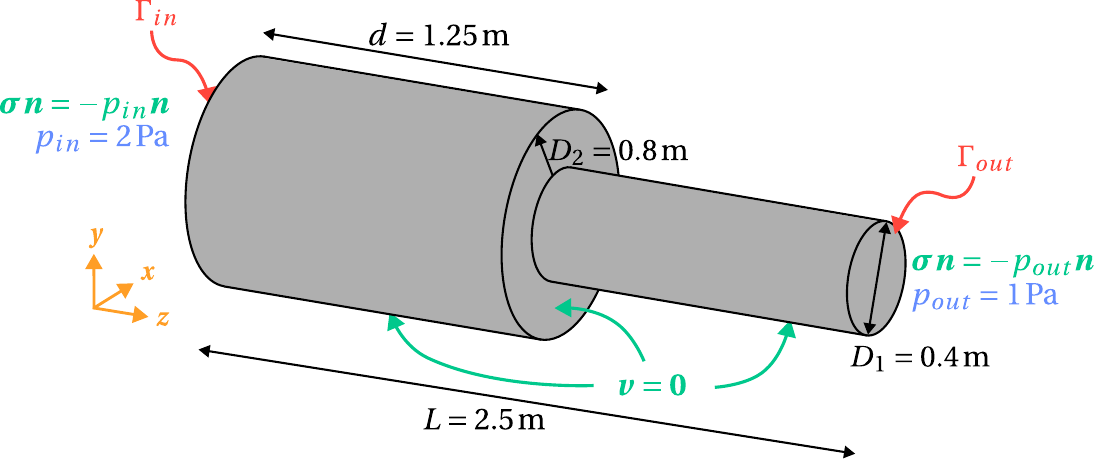}
    \caption{Scheme of the studied problem}
    \label{fig:scheme_reducer}
\end{figure}

The motivation of this section is to demonstrate that the solver also works for 3D problems. We propose to study a reducer composed of two coaxial channels with different diameters. Like in the two previous examples, the reducer is submitted to a pressure $p_{in} = \qty{2}{\pascal}$ on the inflow boundary and $p_{out} = \qty{1}{\pascal}$ on the outflow boundary. We also consider a no-slip boundary condition for the velocity on the external surface of the reducer. The boundary conditions and the geometry of the problem are represented in \autoref{fig:scheme_reducer}.

We use a structured mesh constituted of 12,026 elements which leads to 256,851 DoFs for the velocity part and 11,176 DoFs for the density part.

\begin{figure}[ht!]
    \centering
     \begin{subfigure}[b]{0.99\textwidth}
         \centering
         \includegraphics[width=12.5cm]{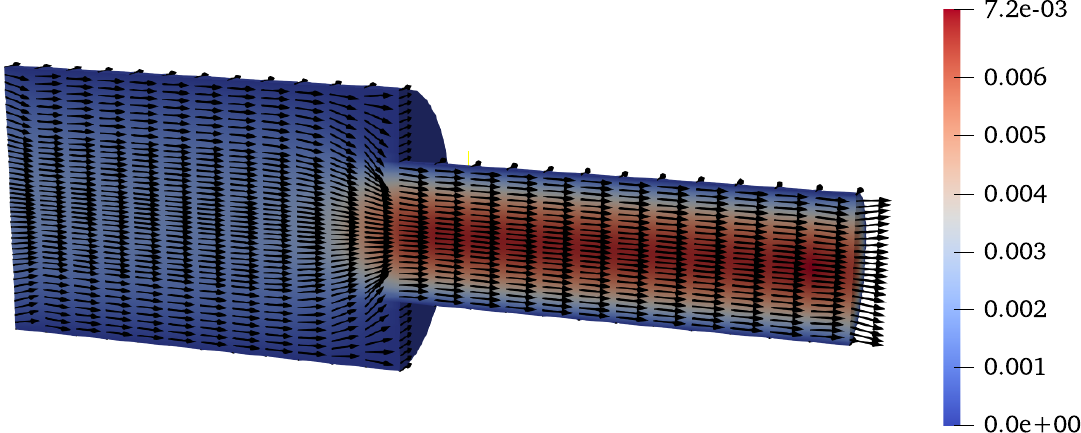}
         \caption{Velocity [m/s]}
         \label{fig:vel_reducer}
     \end{subfigure}
     \hfill
     \begin{subfigure}[b]{0.99\textwidth}
         \centering
         \includegraphics[width=12cm]{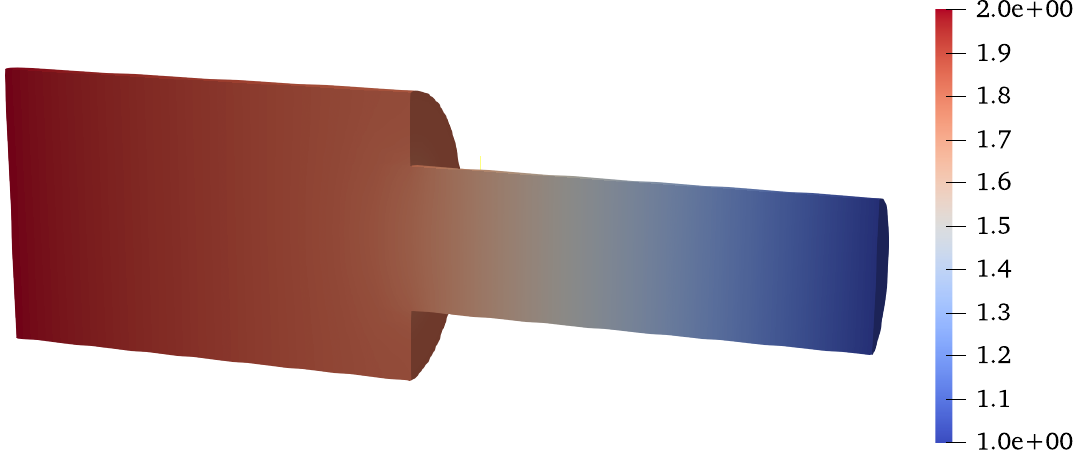}
         \caption{Pressure [Pa]}
         \label{fig:pressure_reducer}
     \end{subfigure}
    \caption{Fields at last time step inside the reducer}
    \label{fig:fields_reducer}
\end{figure}

The results for the velocity and the pressure are depicted in \autoref{fig:fields_reducer}. The pressure regularly decreases along the main axis. For the velocity, we distinctly notice the two areas with different cross sections. When the fluid enters the channel with the smallest cross-section, its velocity increases significantly. Moreover, using streamlines, we observe the viscous effects at the points where the channel diameter changes.

The simulation yields these results in 40 iterations. However, these iterations take longer than the previous two problems as the reducer has a significantly higher number of degrees of freedom. The velocity basis contains 20 modes, while the pressure basis contains 9. Using these bases, we achieve a global relative error of $0.2\%$ for the velocity and $0.2\%$ for the pressure.

Therefore, we have demonstrated that the code works for different 2D and 3D geometries. It is then predictable that the method could also be applied to more complex meshes.

\section{Conclusion}\label{sec:conclusion}

Previous work has used PGD to address flow problems or the LATIN solver to handle the specific case of fluid-structure interaction. However, this work represents the first attempt to apply the LATIN-PGD method to the resolution of fluid problems. The first results demonstrate that this solver effectively computes relevant solutions. We have recovered analytical results for a Poiseuille flow and reproduced results from the literature on a standard benchmark. Moreover, it is also possible to compute problems for 3D geometries. For all these examples, the pressure and velocity have been decomposed into a space-time separated variable form. With a relatively small basis, it is possible to get accurate approximations of a full-order simulation.

However, several assumptions have been made in our models. We restricted ourselves to Newtonian compressible laminar flows. Considering the nonlinear convective term would require including a classical Newton-Raphson scheme at the local stage. We could also tackle the classical case of incompressible flows by drawing inspiration from velocity-pressure decoupling strategies, such as those presented in \cite{wang_2018,decaria_2020}.

More importantly, this work represents just the first step in the field of fluid problems and offers numerous possibilities for future research. As a recall, the LATIN-PGD solver has proven relevant for integrating complex material constitutive laws in solid mechanics. Thus, this solver could be relevant for tackling complex non-Newtonian flow problems. For such problems, the viscosity is no longer constant, and models of varying complexity have been proposed in the literature (\cite{ikoku_1979,roussel_2006,junker_2025}). The solver would also be of interest for parametrised problems, either on the geometry, the material or the boundary conditions.

\section*{Declarations}
\begin{itemize}
    \item Availability of data and materials \\
    No data or material was used for the research described in the article.
    \item Competing interests \\
    The authors declare that they have no known competing financial interests or personal relationships that could have appeared to influence the work reported in this paper.
    \item Authors' contributions \\
    E. F.: Writing – original draft, Software, Methodology, Investigation, Conceptualization. P.-A. B.: Review \& editing, Supervision, Methodology. F. L.: Review \& editing, Supervision, Methodology. D. N.: Review \& editing, Supervision, Methodology. P.J.: Review \& editing, Supervision, Methodology, Investigation, Conceptualization.
    \item Acknowledgements \\
    This work was performed using HPC resources from the ``Mésocentre'' computing center of CentraleSupélec and \'Ecole Normale Supérieure Paris-Saclay supported by CNRS and Région \^Ile-de-France \url{https://mesocentre.universite-paris-saclay.fr/}. The authors gratefully acknowledge the support from the German Research Foundation (DFG) within the International Research Training Group 2657 (IRTG 2657) entitled ``Computational Mechanics Techniques in High Dimensions'' (Grant No. 433082294).
\end{itemize}

\bibliographystyle{elsarticle-num} 
\bibliography{biblio}



\end{document}